\documentclass{amsart}
\usepackage{amsmath,amsthm,amssymb,latexsym,amscd}
\theoremstyle{plain}
\newtheorem{theorem}{Theorem}

\theoremstyle{definition}

\theoremstyle{remark}
\newtheorem{remark}{Remark}

\theoremstyle{remark}

\newcommand{\lat}{\mathbb{L}\mathrm{at}^n_r(K)}

\numberwithin{equation}{section}
\begin{document}
\title[The subregular variety of a lattice variety]{Gorenstein and normal properties of the subregular variety of a variety of special lattices over Witt vectors}
\author{Akira Sano}
\address{Center for Advanced Studies in Mathematics\\
         Ben-Gurion University of the Negev\\
         Be'er-Sheva 84105 Israel}
\email{sano@math.bgu.ac.il}
\thanks{The author was supported by the Center for Advanced Studies in Mathematics, Ben-Gurion University of the Negev, Be'er-Sheva, Israel}
\keywords{lattice varieties, Gorenstein rings, complete intersection rings, Witt vectors}
\subjclass[2000]{14M15, 20G25, 20G99}
\begin{abstract}
We recall the basic geometric properties of the projective variety $\lat$ parametrizing a family of special lattices over Witt vectors proved in \cite{Haboush1}.  It is a `Schubert variety' in the analog of the affine Grassmannian for $G=\mathrm{SL}(n)$ over a complete local field of mixed characteristic, and is normal and Gorenstein \cite{Sano, HS}.  In this paper, we prove that there exists a set of subvarieties of $\lat$ that are also normal and Gorenstein.  The set contains the subregular variety, that is, the complement of the smooth locus, of $\lat$.
\end{abstract}
\maketitle
\section{Introduction}
Let $p$ be a prime, $k$ the algebraic closure of the finite field $\mathbb{F}_p$, $\mathcal O=\mathfrak{W}(k)$ the $k$-points of the Witt vectors, and $K$ the fraction field of $\mathcal O$.  Let $F=\mathcal O^n$ be the standard $\mathcal O$-free submodule of $K^n$ with a basis $\{e_1,\dots,e_n\}$ which is also a $K$-basis.  A \emph{lattice} $L\subset K^n$ is a free maximal rank $\mathcal O$-submodule, it is \emph{special} with respect to $F$ if $\bigwedge^n L=\bigwedge^n F$, and it is of \emph{height} at most $r>0$ if $L\subset p^{-r}F$.  We denote by $\lat$ the set consisting of all special lattices of height at most $r$ with respect to $F$.

Let $G_0=\mathrm{SL}_n(\mathcal O)$ and $B_0$ be the Iwahori subgroup of $G_0$ whose subdiagonal entries are divisible by $p$.  From various parts of \cite{Haboush1}, we have
\begin{theorem}\label{T:lattice}\mbox{}
\begin{enumerate}
\item The set $\lat$ is a projective $k$-variety.
\item For every dominant cocharacter $\gamma=\mathrm{diag}(p^{r_1},\dots,p^{r_n})\in\Gamma^+$, that is, $(n-1)r\geq r_1\geq\dots\geq r_n\geq -r$ and $\sum_{i=1}^n r_i=0$, the $G_0$-orbit closures $\overline{G_0\gamma F}$ (the proalgebraic action defined up to a suitable Frobenius cover of $G_0$) have even dimensions equal to $2\sum_{i=2}^n (1-i)r_i$.\item The smooth locus $\lat_{\rm{reg}}$ is exactly the orbit $G_0\mu_rF$, where $\mu_r$ is the diagonal matrix with respect to the basis $\{e_1,\dots,e_n\}$ with the diagonal entries $p^{(n-1)r},p^{-r},\dots,p^{-r}$.  In particular,  $\dim\lat=\dim G_0 \mu_r F = (n-1)nr$.
\item The orbit $B_0 \mu_r F$ is isomorphic to the affine space $\mathbb{A}_k^{(n-1)nr}$.
\item The complement of $G_0 \mu_r F$ is of codimension two.
\end{enumerate}
\end{theorem}
In \cite{Sano}, the author proved
\begin{theorem}\label{T:mt}
The variety $\lat$ is normal and Frobenius split.
\end{theorem}
In the recent paper \cite{HS}, we presented a shorter and more canonical proof of the normality, and observed that the basic construction in the proof of Theorem~\ref{T:mt} carries the result further; namely,
\begin{theorem}\label{T:gore}
The variety $\lat$ is normal and Gorenstein.
\end{theorem}
In this paper, we will study the subregular variety of $\lat$ as a particular member in a set of subvarieties in $\lat$.  In Theorem~\ref{T:main}, we shall prove the same geometric properties for them as those in Theorem~\ref{T:gore}.  The proof of Theorem~\ref{T:main} as well as that of the corresponding results for $\lat$ in \cite{Sano, HS} are both based on showing the complete intersection property of a certain covering variety.  This approach is similar to the proof of the normality of the nilpotent variety in \cite{KP}.  The geometry of $\lat$ is similar to that of the nilpotent variety, with Gorenstein and normal properties also holding true for the subregular nilpotent variety, consisting of irregular nilpotent elements, of a semisimple Lie algebra proved in e.g., \cite{Broer, KLT}.

Over the complete local field $K$ of equal characteristic, more specifically for $K=k((t))$, it is known that all Schubert varieties in affine (partial) flag varieties for $G(K)$ with $G$ reductive, are normal and Cohen-Macaulay with rational singularities (see, e.g. \cite{Faltings}).  In particular, Cohen-Macaulay property is proved using various techniques including Frobenius splitting, induction on dimension and cohomology vanishing.  Over the complete local field of mixed characteristic, more specifically for $K$ the fraction field of Witt vectors over $k$, the same geometric properties hold for the affine-Grassmannian-like variety with $G=\mathrm{SL}_n$, denoted by $\varinjlim_{r}\lat$ (which is defined over $k=\overline{\mathbb{F}}_p$ and not over $\mathbb{Z}$), using Frobenius splitting and related methods (\cite{Sano2}).  The proofs of Theorems~\ref{T:mt}, \ref{T:gore} and \ref{T:main} of normal and Gorenstein properties of \emph{certain} Schubert varieties in $\varinjlim_{r}\lat$ require neither methods such as induction on the dimension nor cohomology vanishings, while they give a stronger property than being Cohen-Macaulay.  The methods we use,  also apply to the same set of Schubert varieties in the equal characteristic analog.  They provide a class of normal and Gorenstein affine Schubert varieties in the affine Grassmannian associated with $G=\mathrm{SL}_n$.
 
An application of the geometry of $\lat$ is given in \cite{HS} where the canonical generator of its Picard group is computed, and the space of its sections is an $\mathrm{SL}_n(\mathcal{O})$-representation in characteristic $p$.  Similar results are expected for the geometry of (proper) Schubert subvarieties of $\lat$.
\section{Geometry of the matrix cover variety: Summary}
We shall review the basic construction and properties of the matrix cover variety $X_r$ and the morphism $\pi:X_r\to\lat$.  The proofs of the statements in this section are provided in \cite{HS, Sano}.  Nonetheless, except for the proof of the smoothness of $\pi$, they will be provided in the proof of Theorem~\ref{T:main} as a special case.

Let $\mathfrak{v}:\mathcal{O}\to\mathcal{O}$ be the Verschiebung morphism, and define, for $F=\mathcal{O}^n$ and $u_i=\sum_{i=0}^{\infty}\xi(u_{ij})^{p^{-i}}{p}^i$, the morphism $\mathfrak{v}_F:F\to F$ by
\begin{equation}\label{E:versh}
\mathfrak{v}_F(u_1,\dots,u_n)^T :=(\mathfrak{v}(u_1),\dots,\mathfrak{v}(u_n))^T.
\end{equation}
Here, $(u_1,\dots,u_n)^T$ denotes the transpose of the row vector $(u_1,\dots,u_n)$.  As in \cite{Haboush1, HS}, let $\lat$ be identified, via $\mathfrak{v}_F^r$, with $\mathbb{L}(0, nr; (n-1)nr)$, which is the set of lattices $L$ such that $p^{nr}F\subset L\subset F$ with $\dim_k L/p^{nr}F=(n-1)nr$.  It contains the lattice $\mathbb{L}(nr,0,\dots,0)\subset F$, spanned by $\{p^{nr}e_1,e_2,\dots,e_n\}$ which is identified with the lattice $\mu_r F \in \lat$ via $\mathfrak{v}_F^r$.  In what follows, we shall write $\mu_r=\mathrm{diag}(p^{nr},1,\dots,1)$.

We denote by $\overline{\mathcal{O}}:=\mathcal {O}/p^{nr+1}\mathcal{O}$ the finite Witt vectors of length $nr+1$, and shall write the elements of $\overline{\mathcal{O}}$ as $\sum_{i=0}^{nr}\xi(a_i)^{p^{-i}}p^i$, where $\xi:k^{\times}\to\mathcal{O}^{\times}$ is the system of multiplicative representatives, and $p^i$ the residue class of $p^i$ modulo $p^{nr+1}\mathcal{O}$.  Denote by $G:=\mathrm{GL}_n(\overline{\mathcal O})$ a finite-dimensional algebraic group over $k$ (while $G_0$ still denotes $\mathrm{SL}_n(\mathcal{O})$), and by $B \subset G$ the image of the standard Iwahori subgroup of $\mathrm{GL}_n(\mathcal{O})$ under the map induced by the canonical map $\mathcal{O}\to\overline{\mathcal{O}}$.

Let $\Gamma(nr+1)=\{g\in \mathrm{GL}_n(\mathcal{O})\mid g\equiv I\pmod{p^{nr+1}\mathcal{O}}\}$ be a congruence subgroup.  We know $\Gamma(nr+1)$ acts trivially on $\lat$, hence $G\cong \mathrm{GL}_n(\mathcal{O})/\Gamma(nr+1)$ acts algebraically on $\lat$.  This way, we reduce the infinite-dimensional nature of our problems to the ordinary algebraic geometry of finite-dimensional (non-reductive) algebraic group $G$ over $k$.

Let $\operatorname{Mat}_n(\overline{\mathcal{O}})$ be the affine $k$-space of dimension $n^2 (nr+1)$, consisting of $n\times n$-matrices with entries in $\overline{\mathcal{O}}$.   For $A=(a_{i,j})\in \operatorname{Mat}_n(\overline{\mathcal{O}})$, we denote by $S$ the set of variables $\{x_{i,j,s} \mid 1\leq i,j \leq n, 0\leq s \leq nr\}$ for the coordinates of $a_{i,j,s}$ where $a_{i,j}=\sum_{s=0}^{nr}\xi(a_{i,j,s})^{p^{-s}}p^s$.  Let $\det:\operatorname{Mat}_n(\overline{\mathcal{O}})\to\overline{\mathcal{O}}$ be the determinant morphism, and let us define 
\begin{equation}\label{E:X_r}
X_r:=\{A\in \operatorname{Mat}_n(\overline{\mathcal O})\mid \det A=p^{nr}u, \mbox{ for some }u\in\overline{\mathcal{O}}^{\times}\}.
\end{equation}
Then $A\in X_r$ if and only if there exist polynomials $d_0,d_1,\dots,d_{nr}\in k[S]$ such that
\begin{equation}\label{E:determinant}
\det A=\sum_{i=0}^{nr}\xi(d_i(A))^{p^{-i}}p^i=p^{nr}u=p^{nr}\xi(a_0^{p^{nr}})^{p^{-nr}}
\end{equation}
for some $u=\sum_{i=0}^{nr}\xi(a_i)^{p^{-i}}p^i\in\overline{\mathcal O}^{\times}$.  We see from (\ref{E:determinant}) that the ideal $I(X_r)$ of the localization $k[S]_{d_{nr}}$ is defined by the simultaneous vanishings of the $nr$ polynomials 
\[
d_0, d_1, \dots, d_{nr-1}.
\]
Thus, $X_r$ is a closed subvariety of the open set $\operatorname{Spec}(k[S]_{d_{nr}})\subset \operatorname{Mat}_n(\overline{\mathcal O})$.

Also we recall that,
\begin{equation}\label{E:closure}
X_r=\overline{G \mu_r G},
\end{equation}
by showing a series of inclusions,
\begin{equation}
\overline{G \mu_r G}\subset X_r
\subset \bigcup_{\begin{subarray}{c}p^r I \leq\eta\leq \mu_r\\ \eta\in\Gamma^+ \end{subarray}} G \eta G
\subset\overline{G \mu_r G }.
\end{equation}
The proof of (\ref{E:closure}), which is in \cite{Sano}, will appear in the proof of Theorem~\ref{T:main} as a special case.  Using this, we can compute the codimension of $X_r$ to be $nr$ in $\operatorname{Mat}_n(\overline{\mathcal{O}})$ by considering the orbit structure of $X_r$.   Hence $X_r$ is a complete intersection in $\operatorname{Mat}_n(\overline{\mathcal{O}})$ and in particular Gorenstein.  Further, one can show that $X_r$ is nonsingular in codimension one.  Together with Serre criterion for normality, we conclude that $X_r$ is normal.

The morphism
\[
\pi:X_r\to \lat
\]
defined by $A\mapsto AF$ is surjective and can be shown to be smooth (cf. \cite[Proposition~1]{HS}.  Note that in \cite{HS}, $X_r$ is denoted by $\mathcal{B}_{nr}(F)$).  Therefore by a standard theorem in commutative algebra, $\lat$ is normal and Gorenstein (cf. \cite[Theorem~32.2]{Matsumura}).  The fiber of $\pi$ at every point is a $G$-orbit.  By another well-known theorem (cf. \cite[Proposition~6.22]{Borel}), the normality of $\lat$ implies that it is the geometric quotient 
\begin{equation}\label{E:quo}
\lat \simeq X_r/G.
\end{equation}
\begin{remark}\label{D:identify}
We can identify $\overline{B_0\eta F}\subset\lat$ with $\overline{B \eta G}/G$ by (\ref{E:quo}), and similarly $\overline{G_0\eta F}\subset\lat$ with $\overline{G \eta G}/G$, so that $\pi^{-1}(\overline{B_0\eta F})=\overline{B \eta G}$ as well as $\pi^{-1}(\overline{G_0\eta F})=\overline{G \eta G}$.  Since $\pi$ is smooth everywhere, so are their restrictions to these inverse images.
\end{remark}
\section{The main theorem}
We denote by $\lat^{\mathrm{sr}}$, the subregular variety of $\lat$.   It is the complement of the smooth locus of $\lat$.  By Theorem~\ref{T:lattice}, it is precisely given by the $G_0$-orbit closure $\overline{G_0\gamma^{\mathrm{sr}} F}$ where $\gamma^{\mathrm{sr}}$ is identified, via the morphism (\ref{E:versh}), with $\mathrm{diag}(p^{nr-1},p,1,\dots,1)$.  More generally, for all $i=0,1,\dots,\ell$ with $\ell\leq nr/2$, we define $\gamma^{\mathrm{sr}}_i=\mathrm{diag}(p^{nr-i},p^i,1,\dots,1)$, and $\lat^{\mathrm{sr}}_i:=\overline{G_0 \gamma^{\mathrm{sr}}_i F}$.  When $i=0$, it is equal to $\lat$ and when $i=1$, the subregular variety $\lat^{\mathrm{sr}}$.

Our main theorem of this paper is
\begin{theorem}\label{T:main}
For all $i=0,1,\dots,\ell$ such that $\ell\leq nr/2$, the subvarieties $\lat^{\mathrm{sr}}_i$ of $\lat$ are normal and Gorenstein.
\end{theorem}
\begin{proof}
As in the case of $\lat$, first we shall define the matrix cover varieties to $\lat^{\mathrm{sr}}_i$, and prove that they are normal and Gorenstein.  Then by the smoothness of the covering morphism, these properties in question are passed to those of the target varieties, hence proving the theorem.

Let $A_{i,j}$ denote the minor matrix of $A=(a_{i,j})\in \operatorname{Mat}_n(\overline{\mathcal{O}})$ formed by removing its $i$-th row and the $j$-th column.  Let $b,c:\operatorname{Mat}_n(\overline{\mathcal{O}})\to \overline{\mathcal{O}}$ be the functions defined by
\begin{align*}
b(A)&:=a_{1,1}=\sum_{j=0}^{nr}\xi(b_j(A))^{p^{-j}}p^j,\\
c(A)&:=\det(A_{1,1})=\sum_{j=0}^{nr}\xi(c_{j}(A))^{p^{-j}}{p}^j,
\end{align*}
where $\{b_{j},c_{j} \mid j=0,1,\dots,nr\}$ is the set of polynomials in $k[S]_{d_{nr}}$.  Let us denote the valuation map by $v_p:\overline{\mathcal{O}}\to\mathbb{Z}$.

For $i=0,1,\dots,\ell$, we define the subvarieties $X_r^{\mathrm{sr}}(i)$ of $X_r$ as:
\begin{equation}
X_r^{\mathrm{sr}}(i):=\{A=(a_{i,j})\in X_r \mid v_p(c(A))\geq i, \text{ and } v_p(b(A))\leq nr-i \}.
\end{equation}
Then $X_r^{\mathrm{sr}}(i)$ is defined by the ideal
\begin{equation}
\langle d_0,d_1,\dots,d_{nr-1};b_{nr-i+1},\dots,b_{nr},c_0,\dots,c_{i-1} \rangle
\end{equation}
which is generated by $nr+2i$ polynomials in the regular local ring $k[S]_{d_{nr}b_{nr-i}c_{i}}$.  When $i=0$, $X_r^{\mathrm{sr}}(0)=X_r$.

Our claim is that each $X_r^{\mathrm{sr}}(i)$ is a complete intersection in $\operatorname{Mat}_n(\overline{\mathcal{O}})$ that maps smoothly onto $\lat^{\mathrm{sr}}_i$. We do so by computing its dimension by showing that it is an orbit closure, and showing that the codimension in $\operatorname{Mat}_n(\overline{\mathcal{O}})$ matches the number of generators of its ideal in $k[S]_{d_{nr}b_{nr-i}c_{i}}$.

It suffices, by Remark~1, to show that 
\begin{equation}\label{E:set}
X_r^{\mathrm{sr}}(i)=\overline{G\gamma^{\mathrm{sr}}_i G}.
\end{equation}
To check (\ref{E:set}), let $P$ be the subgroup of $G$ whose elements are the matrices of the form:
\[
\begin{pmatrix}
* & *\\
0 & \mathrm{GL}_{n-1}(\overline{\mathcal{O}})
\end{pmatrix},
\]
and $P^-$ be the subgroup whose elements are the matrices of the form:
\[
\begin{pmatrix}
* & 0\\
* & \mathrm{GL}_{n-1}(\overline{\mathcal{O}})
\end{pmatrix}.
\]
Clearly, $\gamma^{\mathrm{sr}}_i\in X_r^{\mathrm{sr}}(i)$.  We have the action of the group $P\times P^-$ on $X_r$ restricted from the action of $G\times G$.  If
\[
X=\begin{pmatrix}x_{1,1}&X_{1,2}\\0&X_{2,2}\end{pmatrix}\in P,\quad Y^{-1}=\begin{pmatrix}y_{1,1}&0\\Y_{2,1}&Y_{2,2}\end{pmatrix}\in P^-,
\]
then
\[
c(X\gamma^{\mathrm{sr}}_i Y^{-1})=\det(X_{2,2}(\begin{smallmatrix}p^i&0\\0&I\end{smallmatrix})Y_{2,2}).
\]
Since $X_{2,2},Y_{2,2}\in\mathrm{GL}_{n-1}(\overline{\mathcal{O}})$,  we have $v_p(c(X\gamma^{\mathrm{sr}}_i Y^{-1}))\geq i$.

Also, we have
\[
b(X\gamma^{\mathrm{sr}}_i Y^{-1})=x_{1,1}p^{nr-i}y_{1,1}+X_{1,2}(\begin{smallmatrix}p^i&0\\0&I\end{smallmatrix})Y_{2,1}.
\]
To show that $v_p(b(X\gamma^{\mathrm{sr}}_i Y^{-1}))\leq nr-i$, first, we have an exercise to check:
\begin{remark}[Exercise]
If $v$ is a valuation and $v(a)\neq v(b)$, then
\[
v(a+b)=\operatorname{min}\{v(a),v(b)\}.
\]

The inequality $\geq$ is by the definition of a valuation.  For the other direction, assuming $v(a)<v(b)$ (with the symmetric argument for $v(a)> v(b)$), we have
\begin{multline}\label{E:ineq}
\operatorname{min}\{v(a),v(b)\}=v(a)=v(a+b-b)\geq \operatorname{min}\{v(a+b),v(-b)\}\\=\operatorname{min}\{v(a+b),v(b)\}=v(a+b)
\end{multline}
since otherwise, (\ref{E:ineq}) would give $v(a)\geq v(b)$, a contradiction.
\end{remark}

Therefore, for $v_p(X_{1,2}(\begin{smallmatrix}p^i&0\\0&I\end{smallmatrix})Y_{2,1}) \neq nr-i$,
\begin{multline}
v_p(b(X\gamma_i^{\mathrm{sr}}Y^{-1})) =\operatorname{min} \{ nr-i, v_p(X_{1,2}(\begin{smallmatrix}p^i&0\\0&I\end{smallmatrix})Y_{2,1}) \}\\
=
\begin{cases}
v_p(X_{1,2}(\begin{smallmatrix}p^i&0\\0&I\end{smallmatrix})Y_{2,1})\lvertneqq nr-i, &\text{if $v_p(X_{1,2}(\begin{smallmatrix}p^i&0\\0&I\end{smallmatrix})Y_{2,1})\lvertneqq nr-i$;}\\
nr-i, &\text{if $v_p(X_{1,2} (\begin{smallmatrix}p^i&0\\0&I\end{smallmatrix}) Y_{2,1})\gvertneqq nr-i$.}
\end{cases}
\end{multline}
If $v_p(X_{1,2}(\begin{smallmatrix}p^i&0\\0&I\end{smallmatrix}) Y_{2,1})=nr-i$, then
\begin{equation}
b(X\gamma^{\mathrm{sr}}_i Y^{-1})=p^{nr-i}(x_{1,1}y_{1,1}+z)
\end{equation}
for some $z\in\overline{\mathcal{O}}$, hence $v_p(b(X\gamma^{\mathrm{sr}}_i Y^{-1}))=nr-i$.  Thus $P \gamma_i^{\mathrm{sr}} P^{-}\subset X_r^{\mathrm{sr}}(i)$, and since $X_r^{\mathrm{sr}}(i)$ is a closed subvariety of $U=\operatorname{Spec}(k[X_r]_{d_{nr}b_{nr-i}c_{i}})$,
\begin{equation}\label{E:con1}
\overline{P \gamma^{\mathrm{sr}}_i P^-}\subset X_r^{\mathrm{sr}}(i),
\end{equation}
with the closure taken in $U$.

For each $a=i,\dots,(n-1)r$, let
\[
\Gamma_a:=\{ \mathrm{diag}(p^{nr-a},p^{r_2},\dots,p^{r_n})\in \Gamma^+\mid \sum_{j=2}^n r_j=a\}.
\]
Then $\cup_{a=i}^{(n-1)r}\Gamma_a\subset X_r^{\mathrm{sr}}(i)$.  We identify the nonzero part of $(n-1)$-tuples $(r_2,\dots,r_n)$ with a partition of $a$.  We will show that $X_r^{\mathrm{sr}}(i)$ is contained in the union of $G\times G$-orbits indexed by $\cup_{a=i}^{(n-1)r}\Gamma_a$.

If $A=(a_{i,j})\in X_r^{\mathrm{sr}}(i)$, we locate a non-zero entry $a_{i,j}=p^{r_n}u$ with the smallest value $r_n$ with $nr\geq r_n\geq 0$ and $u$ a unit.  Locate the largest row index $i'\geq i$ in the column with $a_{i',j}=p^{r_n}u'$ where $r_n$ is achieved and $u'$ a unit.  For every non-zero entry below it in that column, the values are higher than $r_n$.  We multiply on the left by a product of elementary matrices in $G$ to clear the entries below $a_{i',j}=p^{r_n}u'$.  For the non-zero entries above $a_{i',j}$, we multiply on the left by elementary matrices in $G$ to clear them.  We may do these operations for every column, that is, leave aside the $j$-th column in which we have $a_{i',j}$ in the $i'$-th row and $0$ elsewhere,
and repeat the same argument.  We obtain the entry with the value $r_{n-1}$ so that $nr\geq r_{n-1}\geq r_n\geq 0$, and so on.  The result is $gA=m$ where $g\in G$ is the product of elementary matrices in $G$, and $m$ is a monomial matrix with non-zero entries $p^{r_i}u_i$'s with $u_i$ units in some permutation.  We permute the rows and columns of $m$ by some permutation matrices $\sigma_1$ on the left and by $\sigma_2$ on the right of $m$ so that
\[
\sigma_1 g A \sigma_2=\mathrm{diag}(p^{r_1}u_1,\dots,p^{r_n}u_n),
\]
or equivalently
\[
A=(g^{-1}\sigma_1^{-1})\gamma(\mathrm{diag}(u_1,\dots,u_n)\sigma_2^{-1})\in
G\gamma G,
\]
where $\gamma=\mathrm{diag}(p^{r_1},\dots,p^{r_n})$ has the properties
\[
nr\geq r_1\geq\cdots\geq r_n\geq 0,\quad \sum_{j=1}^{n}r_j=nr,\quad \sum_{j=2}^{n}r_j=a
\]
with $i\leq a\leq (n-1)r$.  Hence
\begin{equation}\label{E:con2}
X_r ^{\mathrm{sr}}(i)\subset \bigcup_{a=i}^{(n-1)r}\bigcup_{\gamma_a\in\Gamma_a}G\gamma_a G.
\end{equation}

Next, we shall prove the containment
\begin{equation}\label{E:con3}
\bigcup_{a=i}^{(n-1)r}\bigcup_{\gamma_a\in\Gamma_a}G \gamma_a G \subset \overline{G\gamma^{\mathrm{sr}}_i G}.
\end{equation}
Suppose that $A\in G\eta G$ for some $\eta\in\Gamma_a$, with $i\leq a\leq (n-1)r$.  We want to show by a degeneration argument, $A\in\overline{G\gamma^{\mathrm{sr}}_i G}$.  Since the last closure is $G\times G$-stable, it suffices to show that for every pair $\eta\in\Gamma_a,\eta'\in\Gamma_{a'}$, with $i\leq a'\leq a\leq (n-1)r$ such that 
\begin{align*}
\eta&=\mathrm{diag}(p^{r_1},\dots,p^{r_j},\dots,p^{r_n}),\\
\eta'&=\eta\,\mathrm{diag}(p^b,1\dots,1,p^{-b},1,\dots,1)\\
&=\mathrm{diag}(p^{r_1+b},\dots,p^{r_{j-1}},p^{r_j -b},p^{r_{j+1}},\dots,p^{r_n})
\end{align*}
with $r_{i-1}\geq r_i -b\geq r_{i+1}$, we have
\[
\eta\in\overline{G\eta'G}.
\]

In order to prove this, it now suffices to show that there is a family of matrices parametrized by $k^{\times}$
\begin{equation}
\mathcal{F}:=\{\eta(t)=\begin{pmatrix}p^{r_1}& & & & \\ &\ddots& & & \\p^{r_j -b}\xi(t)& &p^{r_j}& & \\ & & &\ddots& \\ & & & &p^{r_n}\end{pmatrix}\mid t\in k^{\times}, 0\leq b\leq r_{j}\}
\end{equation}
so that $\mathcal{F}\subset G\eta'G$.  Then we would have 
\[
\eta\in\overline{G\eta'G}.
\]
To construct such $\mathcal{F}$, we shall show that for every $t\in k^{\times}$,\begin{equation}\label{E:defo}
\eta(t)=x\eta'y^{-1}
\end{equation}
for some $x,y\in G$.  We will show that $x,y$ are products of elementary unipotent matrices in $G$.
To simplify the matrix notations, first we consider $\eta(t)=(a_{i,j})$ and construct a $2\times 2$-matrix $\begin{pmatrix}a_{1,1}&a_{1,j}\\a_{j,1}&a_{j,j}\end{pmatrix}$.  Then we can reduce the checking of (\ref{E:defo}) to that of the corresponding $2\times 2$-matrix.  We consider, for all $t\in k^{\times}$, the problem of decomposing the matrix
\[
\begin{pmatrix}
p^{r_1}&0\\p^{r_j -b}\xi(t)&p^{r_j}
\end{pmatrix}
\in \operatorname{Mat}_2(\overline{\mathcal{O}})
\]
which corresponds to $\eta(t)$.  We perform elementary row and column operations with matrices over $\overline{\mathcal{O}}$, and note that for all $t\neq 0$, we have $\xi(t)^{-1}=\xi(t^{-1})$.  Then we can simply check that by computation
\begin{multline}\label{E:fac}
\begin{pmatrix}p^{r_1}&0\\p^{r_j -b}\xi(t)&p^{r_j}\end{pmatrix}\\
=\begin{pmatrix}1&-p^{r_1+b-{r_j}}(p^b-1)\xi(t)^{-1}\\0&1\end{pmatrix}
\begin{pmatrix}p^{r_1 +b}&0\\0&p^{r_j -b}\end{pmatrix}\\
\begin{pmatrix}1&0\\ \xi(t)&1\end{pmatrix}
\begin{pmatrix}1&(p^b-1)\xi(t)^{-1}\\ 0&1\end{pmatrix}.
\end{multline}
Therefore the second matrix factor on RHS of (\ref{E:fac}) corresponds to $\eta'$ while the rest on RHS corresponds to the elements in $G$.  In particular, we take $\eta'=\gamma^{\mathrm{sr}}_i$ so that $r_1 +b=nr-i$.  This proves the containment (\ref{E:con3}).

To finish the proof of (\ref{E:set}), we need to prove
\begin{equation}\label{E:PxP-}
\overline{P\gamma^{\mathrm{sr}}_i P^-}=\overline{G\gamma^{\mathrm{sr}}_i G}.
\end{equation}
Since the LHS is contained in the RHS, it suffices, by showing that
\begin{equation}\label{E:equal}
\dim P\gamma^{\mathrm{sr}}_i P^-=\dim G\gamma^{\mathrm{sr}}_i G=\dim
\operatorname{Mat}_n(\overline{\mathcal{O}})-(nr+2i).
\end{equation}
Let $x, y\in G$.  Write
\begin{equation}\label{E:matrices}
x=\begin{pmatrix}a&b&c\\ d&e&f\\ g&h&J\end{pmatrix},\quad y=\begin{pmatrix}a'&b'&c'\\ d'&e'&f'\\ g'&h'&J'\end{pmatrix}.
\end{equation}
Then $x\gamma_i^{\mathrm{sr}} y^{-1}=\gamma_i^{\mathrm{sr}}$ if and only if $x\gamma_i^{\mathrm{sr}}=\gamma_i^{\mathrm{sr}} y$, that is,
\[
\begin{pmatrix}a&b&c\\ d&e&f\\ g&h&J\end{pmatrix}\begin{pmatrix}p^{nr-i}& & \\ &p^{i}& \\ & & I\end{pmatrix}=\begin{pmatrix}p^{nr-i}& & \\ &p^{i}& \\ & & I\end{pmatrix}\begin{pmatrix}a'&b'&c'\\ d'&e'&f'\\ g'&h'&J'\end{pmatrix}
\] 
that is,
\[
\begin{pmatrix}p^{nr-i}a&p^{i}b&c\\p^{nr-i}d&p^{i}e&f\\p^{nr-i}g&p^{i}h&J\end{pmatrix}=\begin{pmatrix}p^{nr-i}a'&p^{nr-i}b'&p^{nr-i}c'\\p^{i}d'&p^{i}e'&p^{i}f'\\g'&h'&J'\end{pmatrix}.
\]
It implies that, for $i\leq nr/2$ there exist block matrices over $\overline{\mathcal{O}}$ of obvious sizes $b_0,c_0,f_0,d'_0,g'_0,h'_0$ such that
\[
(x,y)=(\begin{pmatrix}a&p^{nr-2i}b_0&p^{nr-i}c_0\\d&e&p^{i}f_0\\g&h&J\end{pmatrix},\begin{pmatrix}a+p^{i+1}\alpha&b'&c'\\p^{nr-2i}d'_0&e+p^{nr+1-i}\beta &f'\\p^{nr-i}g'_0&p^{i}h'_0&J\end{pmatrix}).
\]

We can solve the matrix equation $x\gamma_i^{\mathrm{sr}}=\gamma_i^{\mathrm{sr}} y$, and it follows that 
\begin{align*}
\begin{split}b_0=b'+p^{i+1}\nu,\quad c_0=c'+p^{i+1}\rho,\quad f_0=f'+p^{nr+1-i}\varepsilon,\\d'_0=d+p^{i+1}\eta,\quad g'_0=g+p^{i+1}\varphi,\quad h'_0=h+p^{nr+1-i}\psi,
\end{split}
\end{align*}
for some block matrices $\nu,\rho,\varepsilon,\eta,\varphi,\psi$.  Then
\begin{multline}
(x,y)=(\begin{pmatrix}a&p^{nr-2i}b'+p^{nr+1-i}\nu&p^{nr-i}c'\\d&e&p^{i}f'\\g&h&J\end{pmatrix},\\ \begin{pmatrix}a+p^{i+1}\alpha&b'&c'\\p^{nr-2i}d+p^{nr+1-i}\eta&e+p^{nr+1-i}\beta&f'\\p^{nr-i}g&p^{i}h&J\end{pmatrix}).
\end{multline}
We add the maximal possible dimensions each entries
\[
a,b',c',d,e,f',g,h,J,\alpha,\beta,\eta\mbox{ and } \nu
\]
can have, so that
\begin{align*}
\dim\mathrm{Stab}_{G\times G}(\gamma_i^{\mathrm{sr}})
&=\dim G + \overbrace{i}^{\nu}+\overbrace{(nr-i)}^{\alpha}+\overbrace{i}^{\eta}+\overbrace{i}^{\beta}\\ &=\dim G+nr+2i.
\end{align*}
It follows that
\begin{align*}
\dim G\gamma^{\mathrm{sr}}_i G &=2\dim G - (\dim G+nr+2i)\\ &=\dim G-(nr+2i)\\
&=\dim \operatorname{Mat}_n(\overline{\mathcal{O}})-(nr+2i).
\end{align*}

We can compute $\dim\mathrm{Stab}_{P\times P^-}(\gamma_i^{\mathrm{sr}})$ by the same manner by letting, in (\ref{E:matrices}), $x\in P, y^{-1}\in P^-$, that is, $b'=c'=d=g=0$ such that $x\gamma_i^{\mathrm{sr}}y^{-1}=\gamma_i^{\mathrm{sr}}$.  The corresponding description of the arbitrary elements of $\mathrm{Stab}_{P\times P^-}(\gamma_i^{\mathrm{sr}})$ is
\[
(x,y)=(\begin{pmatrix}a&p^{nr-i+1}b & 0\\0 & e & p^{i}f'\\0 & h & J\end{pmatrix}, \begin{pmatrix}a+p^{i+1}\alpha & 0 & 0\\p^{nr-i+1}d' & e+p^{nr-i+1}\beta & f'\\ 0 & p^i h & J\end{pmatrix}).
\]
Then
\begin{align*}
\dim\mathrm{Stab}_{P\times P^-}(\gamma_i^{\mathrm{sr}})
=&\overbrace{(nr+1)}^{a}+\overbrace{i}^{b}+\overbrace{(nr+1)}^{e}+\overbrace{(n-2)(nr+1)}^{h}\\
&+\overbrace{(n-2)^2(nr+1)}^{J}+\overbrace{i}^{d'}+\overbrace{(n-2)(nr+1)}^{f'}+\overbrace{(nr-i)}^{\alpha}+\overbrace{i}^{\beta}\\
=&(nr+1)[(n-1)^2 +1] + nr+2i.
\end{align*}
Since $\dim P=\dim P^-=[n^2-(n-1)](nr+1)=(nr+1)(n^2 -n+1)$,
\begin{align*}
\dim P\gamma^{\mathrm{sr}}_i P^- &=2(nr+1)(n^2 -n+1)-[(nr+1)((n-1)^2 +1) + nr+2i]\\
&=(nr+1)n^2 -(nr+2i)\\
&=\dim \operatorname{Mat}_n(\overline{\mathcal{O}})-(nr+2i).
\end{align*}
Because $nr+2i$ is exactly the minimal number of generaters of the ideal $I(X_r^{\mathrm{sr}}(i))$,  $X_r^{\mathrm{sr}}(i)$ is a complete intersection in $\operatorname{Spec}(k[S]_{d_{nr}b_{nr-i}c_{i}})$; in particular, it is Gorenstein.

We know that $\pi:X_r\to\lat$ is smooth, hence so is its restriction to $\pi^{-1}(\overline{G_0 \gamma_i^{\mathrm{sr}} F})=\overline{G \gamma_i^{\mathrm{sr}} G}=X_r^{\mathrm{sr}}(i)$.  The complement of the maximal dense orbit $G_0\gamma^{\mathrm{sr}}_i F$ in $\lat^{\mathrm{sr}}$ is $\overline{G_0\gamma_{i+1}^{\mathrm{sr}}F}$, which is of codimension 2 by Theorem~\ref{T:lattice}, so the inverse image in $X_r^{\mathrm{sr}}(i)$ of the complement is $X_r^{\mathrm{sr}}(i+1)$ which is also of codimension 2.  Hence $X_r^{\mathrm{sr}}(i)$ is regular in codimension one.  It follows that $X_r^{\mathrm{sr}}(i)$ is normal by Serre criterion.  By \cite[Theorem~32.2]{Matsumura}, $\lat^{\mathrm{sr}}_i$ is Gorenstein and normal.  For $i=1$, it is the subregular variety of $\lat$.
\end{proof}
\begin{remark}
Let $B^-$ be the Iwahori subgroup whose super diagonals are in $p\overline{\mathcal{O}}$.  Like (\ref{E:PxP-}), we can show that $\overline{G\gamma_i^{\mathrm{sr}}G}=\overline{B\gamma_i^{\mathrm{sr}} B^{-}}$ below.  Its consequence is
\[
\overline{B \gamma_i^{\mathrm{sr}}G}\subset \overline{G\gamma_i^{\mathrm{sr}}G}=\overline{B\gamma_i^{\mathrm{sr}} B^{-}} \subset \overline{B\gamma_i^{\mathrm{sr}}G},
\]
and it shows that the class of subvarieties we proved to be normal and Gorenstein in Theorem~\ref{T:main} are indeed Schubert varieties $\overline{B\gamma_i^{\mathrm{sr}} G}/G\simeq \overline{B_0\gamma_i^{\mathrm{sr}}F}$ by definition and Remark~\ref{D:identify}.  We follow the same computation as before, computing explicitly the dimension of $\mathrm{Stab}_{B\times B^-}(\gamma_i^{\mathrm{sr}})$.  If
\[
(x,y)=(\begin{pmatrix}a&b&c\\pd&e&f\\pg&ph&J\end{pmatrix},\begin{pmatrix}a'&pb'&pc'\\d'&e'&pf'\\g'&h'&J'\end{pmatrix})\in \mathrm{Stab}_{B\times B^-}(\gamma_i^{\mathrm{sr}}),
\]
then using the same notation as before,
\begin{multline}
(x,y)=(\begin{pmatrix}a&p^{nr-2i+1}b'+p^{nr-i+1}\nu&p^{nr-i+1}c'\\
pd&e&p^{i+1}f'\\pg&ph&J_0\end{pmatrix},\\ \begin{pmatrix}a+p^{i+1}\alpha&pb'&pc'\\p^{nr-2i+1}d+p^{nr-i+1}\eta&e+p^{nr-i+1}\beta&pf'\\p^{nr-i+1}g&p^{i+1}h&J_0\end{pmatrix})
\end{multline}
where $J_0$ is such that $(\begin{smallmatrix}I_2&0\\0&J_0\end{smallmatrix})\in B\cap B^-$.  We compute the dimension of $\mathrm{Stab}_{B\times B^-}(\gamma_i^{\mathrm{sr}})$ to be
\[
\dim B\cap B^{-} + \overbrace{nr-i}^{\alpha}+\overbrace{i}^{\beta}+\overbrace{i}^{\nu}+\overbrace{i}^{\eta} = n^2(nr)+n+nr+2i,
\]
and the dimension of the orbit $B\gamma_i^{\mathrm{sr}}B^-$ as
\[
2[n^2(nr)+n(n+1)/2]-[n^2(nr)+n+nr+2i]=n^2(nr+1)-(nr+2i).
\]
\end{remark}
\bibliographystyle{amsalpha}
\bibliography{subregular}  

\providecommand{\bysame}{\leavevmode\hbox to3em{\hrulefill}\thinspace}
\providecommand{\MR}{\relax\ifhmode\unskip\space\fi MR }
\providecommand{\MRhref}[2]{%
  \href{http://www.ams.org/mathscinet-getitem?mr=#1}{#2}
}
\providecommand{\href}[2]{#2}
\begin{thebibliography}{KLT99}

\bibitem[Bor91]{Borel}
A.~Borel, \emph{Linear {A}lgebraic {G}roups \rm{(Second Enlarged Edition)}},
  Graduate Texts in Mathematics, vol. 126, Springer-Verlag, New York, 1991.

\bibitem[Bro93]{Broer}
A.~Broer, \emph{Line bundles on the cotangent bundle of the flag variety},
  Invent. Math. \textbf{113} (1993), no.~1, 1--20. \MR{MR1223221 (94g:14027)}

\bibitem[Fal03]{Faltings}
G.~Faltings, \emph{Algebraic loop groups and moduli spaces of bundles}, J. Eur.
  Math. Soc. (JEMS) \textbf{5} (2003), no.~1, 41--68. \MR{MR1961134
  (2003k:14011)}

\bibitem[Hab05]{Haboush1}
W.~J. Haboush, \emph{Infinite dimensional algebraic geometry: algebraic
  structures on {$p$}-adic groups and their homogeneous spaces}, Tohoku Math.
  J. (2) \textbf{57} (2005), no.~1, 65--117. \MR{MR2113991}

\bibitem[HS]{HS}
W.J. Haboush and A.~Sano, \emph{The normality of certain varieties of special
  lattices}, in preparation.

\bibitem[KLT99]{KLT}
S.~Kumar, N.~Lauritzen, and J.~F. Thomsen, \emph{Frobenius splitting of
  cotangent bundles of flag varieties}, Invent. Math. \textbf{136} (1999),
  no.~3, 603--621. \MR{2000g:20088}

\bibitem[KP79]{KP}
H.~Kraft and C.~Procesi, \emph{Closures of conjugacy classes of matrices are
  normal}, Invent. Math. \textbf{53} (1979), no.~3, 227--247. \MR{MR549399
  (80m:14037)}

\bibitem[Mat89]{Matsumura}
H.~Matsumura, \emph{Commutative {R}ing {T}heory}, Cambridge Studies in Advanced
  Mathematics, vol.~8, Cambridge University Press, Cambridge, 1989.
  \MR{90i:13001}

\bibitem[San]{Sano2}
A.~Sano, \emph{The normality of orbit closures of a certain variety of special
  lattices}, in preparation.

\bibitem[San04]{Sano}
\bysame, \emph{The geometry of finite lattice varieties over {W}itt vectors},
  Ph.D. thesis, The {U}niversity of {I}llinois at {U}rbana-{C}hampaign, 2004.

\end{thebibliography}
\end{document}